\documentclass[reqno]{amsart}
\usepackage{caption}
\usepackage{subcaption}

\usepackage{mathabx,cite}
\usepackage[utf8]{inputenc} 
\usepackage[english]{babel}
\usepackage{amstext}
\usepackage{euscript}
\usepackage{bbm}
\usepackage{mathrsfs}
\usepackage{enumerate}
\usepackage{graphicx}
\usepackage{amscd}
\usepackage{verbatim}
\usepackage{amssymb}
\usepackage{amsthm}
\usepackage{amsmath}
\usepackage{amstext}
\usepackage{amsfonts}
\usepackage{latexsym}
\usepackage{mathtools} 
\usepackage{enumitem} 
\usepackage{accents} 
\usepackage[foot]{amsaddr} 

 \usepackage[left=1.2in, right=1.2in, top=1.5in, bottom=1.5in]{geometry}

\usepackage[usenames,dvipsnames]{xcolor} 
\definecolor{dkblue}{RGB}{1,31,91} 
\usepackage[colorlinks=true, pdfstartview=FitV, linkcolor=dkblue, citecolor=dkblue, urlcolor=dkblue]{hyperref}

\usepackage{todonotes}

\newcommand{\RR}{\mathbb R}

\newcommand{\ZZ}{\mathbb Z}

\newcommand{\TT}{\mathbb T}

\newcommand{\pax}{\partial_x}

\normalsize
\normalsize
\DeclareMathOperator{\diff}{d\!}

\newcommand{\norm}[1]{\left\lVert#1\right\rVert}

\theoremstyle{definition}
\newtheorem{theorem}{Theorem}

\newtheorem{remark}[theorem]{Remark}

\numberwithin{equation}{section}
\numberwithin{theorem}{section}
\numberwithin{definition}{section}

\begin{document}

\subjclass[2010]{}

\title[On two systems of Burgers type arising in nonlinear wave interactions]{On two systems of Burgers type arising in nonlinear wave interactions}

\author[D. Alonso-Or\'an]{Diego Alonso-Or\'an}
\address{Departamento de An\'{a}lisis Matem\'{a}tico and Instituto de Matem\'aticas y Aplicaciones (IMAULL), Universidad de La Laguna, C/Astrof\'{i}sico Francisco S\'{a}nchez s/n, 38271, La Laguna, Spain. \href{mailto:dalonsoo@ull.edu.es}dalonsoo@ull.edu.es}

\author[R. Granero-Belinch\'on]{Rafael Granero-Belinch\'on}
\address{Departamento  de  Matem\'aticas,  Estad\'istica  y  Computaci\'on,  Universidad  de Cantabria.  Avda.  Los  Castros  s/n,  Santander,  Spain. \href{mailto:rafael.granero@unican.es}rafael.granero@unican.es}

\begin{abstract}
In this note we study the well-posedness of two systems of Burgers type arising in nonlinear wave interactions. The first model describes the interaction of a Burger's bore with the classical Korteweg-de Vries equation while the second exemplify the interaction of weak sound waves and entropy waves with small amplitudes. For the former, we show the local existence and uniqueness of solutions in Sobolev spaces and Wiener-type spaces. For the latter, we provide an elementary proof of finite time singularity. 
\end{abstract}
\thispagestyle{empty}

\maketitle
\tableofcontents

\section{Introduction}
Wave interactions play a crucial role in the dynamics of various fluid systems. More precisely, In the field of hydrodynamics, waves are generated by a range of disturbances in a fluid medium, including wind stress on ocean surfaces, pressure fluctuations, and gravitational influences. These waves can manifest as surface gravity waves, internal waves within stratified fluids, capillary waves, and solitons, each defined by unique restoring forces and patterns of wave behavior. Investigating the interactions between these waves is essential for gaining insights into the transfer of energy, momentum, and mass across varying spatial and temporal scales. A key element of wave interactions in hydrodynamics is nonlinearity. While linear theory can effectively describe basic wave propagation, it often falls short in addressing the complexities of actual wave interactions. Nonlinear phenomena become particularly significant when wave amplitudes are large or when waves travel through highly variable environments. In the present work we are interested in studying two different systems describing the interactions of waves. \medskip

\subsubsection*{The Burgers-KdV model} 
The first model, derived recently in \cite{Darryletal}, examines one-dimensional nonlinear wave-current interactions between the classical Burger's equation \cite{Burgers}
\[ \partial_{t}u+3\partial_{x}u=0, \quad (x,t)\in \mathbb{T}\times [0,T), \]
and the soliton solutions of the Korteweg-de Vries equation \cite{KdV}
\[ \partial_{t}v+6v\partial_{x}v+\gamma \partial_{x}^{3}v=0,  \quad (x,t)\in \mathbb{T}\times [0,T). \]
Using Hamilton’s principle, applied to the sum of the Lagrangians for the Burgers and KdV equations, coupled through the product of the fluid velocity and the wave momentum map, results in the Burgers-KdV system \footnote{The Burgers-KdV system treated in this work is not the one treated in \cite{Jeffrey}. } of partial differential equations that reads
\begin{subequations}\label{BKdV}
\begin{align}
\partial_t u+3u\partial_x u&=-6v^2\partial_x v-\gamma v\partial_x^3 v\\
\partial_t v-6v\partial_x v&=\gamma \partial_x^3 v+\partial_x(uv).
\end{align}
\end{subequations}
with $(x,t)\in \mathbb{T}\times [0,T)$ and $\gamma\in \mathbb{R}$. The right-hand sides of these equations capture the interaction between the two separate equations. System \eqref{BKdV} enjoys a Hamiltonian structure given by the functional 
\[ \mathcal{H}(t)=\int_{\mathbb{R}} \frac{1}{2}u^2+\frac{\gamma}{2}(\partial_x v)^2-v^3dx.\]
 In \cite{Darryletal}, in addition to the derivation of \eqref{BKdV}, the authors show different numerical simulations to \eqref{BKdV} that exhibit a transfer in the momentum between the Burgers bore and the KdV wave creating oscillations and sharp edge for mean-flow velocity $u$. 

\medskip

\subsubsection*{The Majda, Rosales \& Schonbek model } The second system we investigate in this work describes the interaction of weak sound waves and entropy waves with small amplitudes and was studied by the pioneering work of Majda \& Rosales \cite{majda1984resonantly}. In \cite{majda1984resonantly} , the authors present a systematic asymptotic theory for resonantly interacting weakly nonlinear hyperbolic waves in a single space variable and come out with the following asymptotic system (see equations (3.25) and (3.26) there)
\begin{subequations}\label{MR}
\begin{align}
\partial_t \sigma_1+\sigma_1\partial_x \sigma_1&=-\frac{1}{2\pi}\int_{-\pi}^\pi \partial_y\sigma_2((x+y)/2) \sigma_3(y)dy\\
\partial_t \sigma_3-\sigma_3\partial_x \sigma_3&=\frac{1}{2\pi}\int_{-\pi}^\pi \partial_y \sigma_2((x+y)/2)\sigma_1(y)dy.
\end{align}
\end{subequations}
Here $\sigma_1$ and $\sigma_1$ are the unknown zero-mean periodic sound waves while $\sigma_2$ is the profile of the entropy wave. At first level of approximation, such an entropy profile is stationary and then it is assumed to be given. After a change of variables, the resonant interaction of non-dispersive hyperbolic waves in one space dimension can also be written as (see equations (1.7) in Majda, Rosales \& Schonbek \cite{majda1988canonical} with $\partial_y E=K$)
\begin{subequations}\label{MRS}
\begin{align}
\partial_t u+u\partial_x u&=-\frac{1}{2\pi}\int_{-\pi}^\pi \partial_y E(x-y) v(y)dy\\
\partial_t v+v\partial_x v&=\frac{1}{2\pi}\int_{-\pi}^\pi \partial_y E(y-x) u(y)dy.
\end{align}
\end{subequations}
We call this latter form of writing the system \eqref{MR} the (MRS) system. Here $u$ and $v$ are the zero-mean $2\pi-$periodic sound waves while $E$ is the profile of the entropy wave. Similar equations also arise in the case of two magnetosonic waves propagating in opposite directions when reflecting resonantly off an entropy wave (see equations (5.1) in Ali \& Hunter \cite{ali1998wave}) or when studying the resonant interaction of small amplitude sound waves with a large amplitude entropy wave with variations that occur only in small regions (see equations (1.3) in Ali \& Hunter \cite{ali2000resonant}). We would like to emphasize that, for a sawtooth profile
\begin{equation}\label{kernel}
E(z)=\left\{\begin{array}{cc}x+\pi & -\pi<x<0
\\ x-\pi & 0<x<\pi 
\end{array}\right.
\end{equation}
the system \eqref{MRS} reduces to
\begin{subequations}\label{MRSLocal}
\begin{align}
\partial_t u+u\partial_x u&=v\\
\partial_t v+v\partial_x v&=-u.
\end{align}
\end{subequations}
Since the derivation forty years ago, the previous systems have been studied both analytically and numerically by many different authors. For instance, Majda, Rosales \& Schonbek \cite{majda1988canonical} proved the existence of smooth periodic wave trains while Pego \cite{pego1988some} established the existence of time periodic ''crossing wave'' trains. We also mention the work by Shefter \& Rosales \cite{shefter1999quasiperiodic} study a class of, so called, ''non breaking for all times'' solutions. There solutions do not show wave breaking leading to shock formation. Also, recently, Qu \& Xin \cite{qu2017global} proved existence of a unique global entropy weak solution for \eqref{MR}. Futhermore, these authors also establish the finite time singularity for the \eqref{MR} in the case where the entropy wave $\sigma_2\in C^1$ and when the sound waves $\sigma_1$ and $\sigma_3$ are stronger than the entropy wave. Finally, Hunter \& Smothers \cite{hunter2019resonant} derived a quasilinear degenerate Schrodinger equation as an asymptotic model of the same physical situation.

\subsection{Notations and functional spaces}
We will use the following notation throughout the manuscript. We are working on the one dimensional torus $\TT=[-\pi,\pi]$, endowed with periodic boundary conditions. For  $1\leq p<\infty$ we denote by  $L^p(\TT;\RR)$ the standard Lebesgue space of measurable $p$-integrable $\RR$-valued functions and by  $L^\infty(\TT;\RR)$ the space of essentially bounded functions. Particularly, $L^2(\TT;\RR)$ is equipped with the inner product 
	$
	(f,g)_{L^2}=\int_{\TT}f\cdot\overline{g}\,{\rm d}x,
	$
	where $\overline{g}$ denotes the complex conjugate of $g$. The Fourier
	coefficients and the Fourier series of $f(x)\in L^2(\TT;\RR)$ are defined by
	$\widehat{f}(\xi)=\int_{\TT}f(x){\rm e}^{-{\rm i}x \xi}\,{\rm d}x$ and 
	$f(x)=\frac{1}{2\pi}\sum_{k\in{\mathbb Z}}\widehat{f}(k){\rm e}^{{\rm i}x k}$, respectively. For $s\in\ZZ^{+}$, the Sobolev space $W^{s,p}(\TT;\RR)$ is defined as
$$ W^{s,p} (\TT)=\{f\in L^{p}(\TT), \partial_{x}^{s}f\in L^{p}(\TT)\} $$
endowed with the norm 
$$\norm{f}^{p}_{W^{s,p}(\TT)}= \norm{f}^{p}_{L^p}+\norm{\partial_{x}^{s}f}^{p}_{L^p}$$
where $\partial_{x}$ is the classical first order differential operator with respect $x$. Their homogenous counterparts will be denoted by  $\dot{W}^{s,p} (\TT,\RR)$ and $\dot{H}^{s}(\TT,\RR)$ when $p=2$. To simplify notation, we will just write
$$ L^p=L^p(\TT,\RR), \quad H^s=H^s(\TT,\RR),\quad \dot{H}^s=\dot{H}^s(\TT,\RR).$$ 

We also introduce the \textit{Wiener space} of analytic function
\[ \mathbb{A}^{s}_{\mu}=\{ f\in L^{2}: (1+|k|^{s})e^{\mu |k|}\widehat{f}(k)\in L^{1} \} \]
endowed with the norm
\[ \norm{f}_{\mathbb{A}^{s}_{\mu}}=\norm{(1+|k|^{s})e^{\mu |k|}|\widehat{f}(k)}_{L^{1}}=\displaystyle\sum_{k\in\mathbb{Z}} (1+|k|^{s}) e^{\mu |k|}|\widehat{f}(k)|.\]
Similarly as before their homogeneous counterparts will be denoted by $\mathbb{\dot{A}}^{s}_{\mu}$. \medskip

We denote with $C=C(\cdot)$ any positive universal constant that may depend on fixed parameters and controlled quantities. Note also that this constant may vary from line to line. It is also important to remind that in order to light notation the condition almost everywhere (a.e) is not always indicated.

 \subsection{Main results}
 As we mentioned in the introduction, the purpose of this work is twofold: first, we will provide a well-posedness result for the Burgers-KdV system \eqref{BKdV} in Sobolev and Wiener spaces. Next, we will provide  an elementary proof of finite time singularity for the previous system MRS \eqref{MRSLocal}. \medskip
 
 In order to state the first result, we will include the effects of dissipation into equation \eqref{BKdV}, this is we consider 
 \begin{subequations}\label{diss:BKdV}
\begin{align}
\partial_t u+3u\partial_x u&=-6v^2\partial_x v-\gamma v\partial_x^3 v+\nu\partial_{x}^{2}u \label{diss:BKdV:u} \\
\partial_t v-6v\partial_x v&=\gamma \partial_x^3 v+\partial_x(uv) \label{diss:BKdV:v}
\end{align}
\end{subequations}
where $\nu>0$ denotes the viscosity parameter. The next theorem, shows the existence and uniqueness of solutions in Sobolev spaces under certain assumptions on the parameters $\gamma$ and $\nu$, namely 
\begin{theorem}\label{teobkdv1}
Let $\nu,\gamma>0$. For $(u_0,v_0)\in H^3\times H^4$ with zero mean initial data, there exists a time $0<T_{max}$ and a unique solution
\[ (u,v)\in C([0,T_{max}); H^3\times H^4)\]
of the initial value problem of \eqref{diss:BKdV} with $u(x,0)=u_{0}(x), v(x,0)=v_{0}(x), \ x\in\mathbb{T}$.
\end{theorem}
\begin{remark}
As the reader will notice in the proof, the dissipation term included in \eqref{diss:BKdV} is essential to close the energy estimates. For $\nu=0$ (and $\gamma\in\mathbb{R}$) it is not clear whether the system is locally well-posed in Sobolev spaces. Furthermore, we believe that the Sobolev regularity exponents $ H^3\times H^4$ is not sharp and that may be slightly  improved using fractional Sobolev spaces and refined commutators. 
\end{remark}

\begin{remark}
Notice that for initial datum $v(x,0)=0$, the evolution equation \eqref{diss:BKdV:v} shows that $v(x,t)\equiv 0$. Hence system \eqref{diss:BKdV} boils down into the classical dissipative Burger's equation where the global existence of solutions is guaranteed for $\nu>0$ as well as  finite time singularities in the inviscid case $\nu=0$.
\end{remark}

As remarked previously, it is not clear whether system \eqref{diss:BKdV:v} is well-posed in Sobolev spaces. However, in the following  result we show that for $\gamma=0, \ \nu=0$ there exists a local in time solution in analytic space $\mathbb{A}_{\mu(t)}$. The precise statement of the result reads
\begin{theorem}\label{teobkdv2}
Let $\nu=\gamma=0$. For $(u_0,v_0)\in \mathbb{A}_{\mu_0}\times \mathbb{A}_{\mu_0}$ with zero mean initial data, there exists a time $0<T_{max}$ and a unique solution
\[ (u,v)\in L^{\infty}((0,T_{max}); \mathbb{A}_{\mu(t)}\times \mathbb{A}_{\mu(t)})\]
of the initial value problem of \eqref{diss:BKdV} with $u(x,0)=u_{0}(x), v(x,0)=v_{0}(x), \ x\in\mathbb{T}$.
\end{theorem}

\begin{remark}
We observe that for $\mu>0$ then the space $\mathbb{A}_{\mu}$ contains analytic functions. This Banach space has been used recently in \cite{GGS20,GS19} to show the local existence of solutions where no dissipation mechanism is available. 
\end{remark}

 To conclude this section, we present the result we provide for the local MRS model 
\begin{subequations}\label{MRSLocal:2}
\begin{align}
\partial_t u+u\partial_x u&=v, \\
\partial_t v+v\partial_x v&=-u.
\end{align}
\end{subequations}
that shows the existence of finite time singularities. More precisely,
\begin{theorem}\label{teomrs}
Let us consider a pair of functions $(u_{0},v_{0})$ satisfying that 
\begin{itemize}
\item $(u_{0},v_{0})$ are odd functions and $
\partial_x u_0(0)-\partial_x v_0(0)=1$,
\item $\partial_x v_0(0)<-1,$
\end{itemize}
Then the solution $(u,v)$ to \eqref{MRSLocal:2}  with initial data $u(x,0)=u_{0}(x), v(x,0)=v_{0}(x)$ is not globally smooth.
\end{theorem}

\begin{remark}
While for the case of $C^1(\mathbb{T})$ kernels $E$ the singularity formation was established recently by Qu \& Xin \cite{qu2017global}, such a finite time singularity for non-smooth kernels is, to the best of the authors knowledge, new. Then, although the singularity result in \cite{qu2017global} does not apply for non-smooth kernels such as \eqref{kernel}, in this paper we show that the result suggested there is correct. In fact, the question of finite time singularity has appeared several times mentioned in the literature. For instance, already in \cite{majda1988canonical} is stated that
\begin{quote}
\textit{In one of the numerical experiments with the local equation from (1.6), we find that smooth initial data with a sufficiently small amplitude never develop shocks throughout a long time interval of integration.}
\end{quote}
\end{remark}


%
%

 \subsection{Outline}
 The rest of the paper is organized as follows.  Section \ref{proof:section:wp} is devoted to show the local-well posedness results for the Burgers-KdV equation. More precisely, in Subsection \ref{proof:section:wp2} and Subsection \ref{proof:section:wp2} we provide the proofs of Theorem \ref{teobkdv1} and Theorem \ref{teobkdv2}, respectively. In Section \ref{proof:section:finite} we demonstrate the finite time singularity of solutions for the local MRS system, stated in Theorem \ref{teomrs}.

\section{Proof of Theorem \ref{teobkdv1} and Theorem \ref{teobkdv2}}\label{proof:section:wp}
In this section we provide the proofs of Theorem \ref{teobkdv1} and Theorem \ref{teobkdv2}.

\subsection{Proof of Theorem \ref{teobkdv1}}\label{proof:section:wp1}
Without loss of generality we assume that $\gamma=\nu=1$ (the proof can be easily adapted to different, positive values of $\nu$ and $\gamma$). The proof follows from the combination of appropriate a priori energy estimates and the use of a suitable approximation procedure using the periodic heat kernel, \cite{MajdaBertozzi}. Thus, we first focus in deriving a priori energy estimates and later comment briefly on the approximation procedure to construct the solution. \medskip

\subsubsection*{A priori energy estimates} Integrating system \eqref{diss:BKdV} in space and using the periodicity we find that
\[  \int_{\mathbb{T}} u(x,t) \ \diff x =\int_{\mathbb{T}} u_{0}(x) \ \diff x, \quad \int_{\mathbb{T}} v(x,t) \ \diff x=\int_{\mathbb{T}} v_{0}(x) \ \diff x, \]
thus, the mean property is conserved in time. It is then sufficient to focus on the higher-order estimates $(u,v)\in \dot{H}^3\times \dot{H}^4$. The evolution of the $\dot{H}^3$ norm is given by
\begin{align*}
\frac{d}{dt}\|\partial_x^3 u\|_{L^2}^2&=-\|\partial_x^4 u\|_{L^2}^2+3\int_{-\pi}^\pi \pax^2(u\pax u)\pax^4 u\ \diff x-6\int_{-\pi}^\pi\pax^3(v^2\pax v)\pax^3 u\ \diff x -\int_{-\pi}^\pi\pax^3(v\pax^3v)\pax^3u\ \diff x \\
&\hspace{0.5cm}= -\|\partial_x^4 u\|_{L^2}^2+I_1+I_2+I_3.
\end{align*}
Using integration by parts and H\"older inequality we obtain that
\begin{equation}
 |I_1|\leq C\|\pax u\|_{L^\infty}\|\pax^3 u\|_{L^2}^2, \quad |I_2| \leq C\|\pax v\|^2_{L^\infty}\|\pax^4 v\|_{L^2}\|\pax^3 u\|_{L^2}  \label{eq:I1}
 \end{equation}

Similarly, we compute the evolution of the $\dot{H}^4$ norm of $v$ 
\begin{align*}
\frac{d}{dt}\|\partial_x^4 v\|_{L^2}^2&= 6\int_{-\pi}^\pi\pax^4(v\pax v)\pax^4 v\ \diff x+\int_{-\pi}^\pi\pax^4v\pax^7v \ \diff x+\int_{-\pi}^\pi\pax^5(uv)\pax^4 v\ \diff x.
 \\
&\hspace{0cm}=J_1+J_2+J_3.
\end{align*}
Once again, the Burgers type term can be bounded by
\begin{equation}\label{eq:J1}
|J_1|\leq C\|\pax v\|_{L^\infty}\|\pax^4 v\|_{L^2}^2
\end{equation}
and $J_2$ vanishes using the periodicity condition since
\[ J_{2}=-\frac{1}{2}\int_{-\pi}^\pi \pax\left(\pax^5v\right)^{2}  \diff x=0.\]
The remaining terms $I_{3}$ and $J_{3}$ must be treated more carefully. We first decompose $J_{3}$ as 
\[ J_{3}=\displaystyle\sum_{j=1}^{6} K_{j},\]
 where
$$
K_1=\int_{-\pi}^\pi \pax^5 u v\pax^4 v \ \diff x, \quad K_2=\int_{-\pi}^\pi u\pax^5 v\pax^4 v \ \diff x, \quad K_3=\int_{-\pi}^\pi \pax^4 u \pax v\pax^4 v \ \diff x,
$$
and
$$
K_4=\int_{-\pi}^\pi \pax^3 u \pax^2 v\pax^4 v \ \diff x, \quad K_5=\int_{-\pi}^\pi \pax^2 u \pax^3 v\pax^4 v \ \diff x, \quad K_6=\int_{-\pi}^\pi \pax u (\pax^4 v)^2 \ \diff x.
$$
We first notice that using H\"older inequality and Gagliardo-Nirenberg interpolation inequality \[ \norm{f}_{L^4}\leq \norm{f}^{\frac{1}{2}}_{L^{\infty}}\norm{f_{x}}^{\frac{1}{2}}_{L^{2}},\]
to estimate $K_{4}$ and $K_{5}$ we infer that
\begin{equation}\label{eq:Ktotal}
|K_{2}+K_{4}+K_{5}+K_{6}|\leq C\left( \|\pax u\|_{L^\infty}\|\pax^4 v\|^{2}_{L^2}+\|\pax^{3} u\|_{L^2}\|\pax^4 v\|_{L^2}\norm{\pax^2 v}_{L^{\infty}}\right).
\end{equation}
On the other hand, recalling the term $I_{3}$ and integrating twice by parts we find that
 \[
I_3=-\int_{-\pi}^\pi \pax^5 u(\pax v\pax^3v+ v\pax^4 v) \diff x,
\] 
and hence combining $I_3+K_1$ we infer that 
$$
I_3+K_1=
 \int_{-\pi}^\pi \pax^4 u\pax^4v\pax v \diff x +\int_{-\pi}^{\pi} \pax^4 u\pax^2v\pax^{3} v \diff x.  $$
Therefore, we have that
\[ |I_{3}+K_{1}|\leq C\|\pax^4 u\|_{L^2}\|\pax^4 v\|_{L^2}\left(\|\pax v\|_{L^\infty}+\|\pax^2 v\|_{L^\infty}\right).\]
Notice that the previous structure in $I_3+K_1$ can be obtained as long as $\gamma>0$. Similarly,
$$
|K_{3}|\leq C\|\pax^4 u\|_{L^2}\|\pax^4 v\|_{L^2}\|\pax v\|_{L^\infty}.
$$
In particular, using $\varepsilon$-Young's inequality, we can find
\begin{equation}\label{eq:Kcom}
|I_3+K_1+K_3|\leq \varepsilon \|\pax^4 u\|_{L^2}^2+ C_{\varepsilon}\|\pax^4 v\|_{L^2}^2\left(\|\pax v\|_{L^\infty}^2+\|\pax^{2} v\|_{L^\infty}^2\right),
\end{equation}
for $\varepsilon>0$ arbitrary.  Therefore, collecting estimates \eqref{eq:I1}-\eqref{eq:Kcom}, using the Sobolev embedding $H^{\frac{1}{2}+\epsilon}\hookrightarrow L^{\infty}$, for $\epsilon>0$ and Young's inequality we conclude that
\begin{align*}
\frac{d}{dt}(\|u\|_{H^3}^2+\|v\|_{H^4}^2)+(1-\varepsilon)\norm{\pax^4 u}_{L^{2}}^{2}&\leq C_{\varepsilon}\left(\|\pax^4 v\|_{L^2}^2+\|\pax^3 u\|_{L^2}^2\right)\left(\|\pax v\|_{L^\infty}+\|\pax^2 v\|_{L^\infty}+\|\pax u\|_{L^\infty}\right) \\
&\leq C_{\varepsilon}(\|u\|_{H^3}^2+\|v\|_{H^4}^2)^2.
\end{align*}
Therefore, denoting by $\mathcal{E}(t)=\|u(t)\|_{H^3}^2+\|v(t)\|_{H^4}^2$ we find that
\[ \mathcal{E}(t)\leq \frac{1}{\mathcal{E}(0)-Ct}, \quad 0<t<T_{\max}, \]
which ensures a local uniform time of existence for the solution 
\[ (u,v)\in L^{\infty}\left([0,T_{\max}); H^{3}\times H^{4}\right).\]

\subsubsection*{The approximation procedure} The local existence of solution will follow now from a standard application of Picard’s theorem to a sequence of approximate problems via the periodic heat kernel, \cite{MajdaBertozzi}. Indeed, these sequence of regularized problems have a unique solution $(u^{\epsilon},v^{\epsilon})=C^{1}([0,T_{\epsilon}]; H^{3}\times H^{4})$. By repeating the same energy estimates, we can take $T_{\max}=T(\mathcal{E}_{0})$ independent of $\epsilon$ and passing to the limit we conclude the existence of at least one solution in
\[(u,v)\in L^{\infty}\left([0,T_{\max}); H^{3}\times H^{4}\right),\]
that we skip for the sake of brevity.

\subsubsection*{Uniqueness} At this level of regularity, the uniqueness of such local strong solution can be easily obtained from a standard contradiction argument. Indeed, assume there exists two soutions $(u_{1},v_{1}),(u_{2},v_{2})$ corresponding to the same initial data $(u_{0},v_{0})$. Denoting by $(\overline{u},\overline{v})=(u_{1}-u{2},v_{1}-v_{2})$ we find that
 \begin{subequations}\label{diss:BKdV:diff}
\begin{align}
\partial_t \overline{u}+\frac{3}{2}\partial_{x}\left(u_{1}^{2}-u_{2}^{2}\right)&=-2\partial_{x}\left(v_{1}^{3}-v_{2}^{3}\right)-\gamma\left( v_{1}\partial_x^3 v_{1}-v_{2}\partial_x^3 v_{2}\right)+\nu\partial_{x}^{2}\overline{u},  \\
\partial_t \overline{v}-3\partial_{x}\left(v_{1}^{2}-v_{2}^{2}\right)&=\gamma \partial_x^3 \overline{v}+\partial_x(u_{1}v_{1})-\partial_{x}(u_{2}v_{2}).
\end{align}
\end{subequations}
Then, multiplying by $(\overline{u},\overline{v})$, integrating by parts and repeating the previous estimates  we have that
\[ \frac{1}{2}\frac{d}{dt}\left(\| \overline{u}\|_{L^2}^2+\|\overline{v} \|_{L^2}^2 \right)+\nu \norm{\partial_{x}\overline{u}}_{L^{2}}^{2}\leq C\beta(t)\left(\| \overline{u}\|_{L^2}^2+\|\overline{v} \|_{L^2}^2 \right)\]
with 
\[\beta(t)=\left(\|\pax v_{1}\|_{L^\infty}+\|\pax^2 v_{1}\|_{L^\infty}+\|\pax u_{1}\|_{L^\infty}+\|\pax v_{2}\|_{L^\infty}+\|\pax^2 v_{2}\|_{L^\infty}+\|\pax u_{2}\|_{L^\infty}\right).\]
Invoking Gr\"onwalls inequality,
\[ \left( \| \overline{u}\|_{L^2}^2+\|\overline{v} \|_{L^2}^2 \right) \leq 2 \left( \| \overline{u}_{0}\|_{L^2}^2+\|\overline{v}_{0} \|_{L^2}^2\right) \mbox{exp}\big(\int_{0}^{t}\beta(s) \ \diff s \big), \]
and hence by means of the Sobolev embedding $H^{\frac{1}{2}+\epsilon}\hookrightarrow L^{\infty}$, the uniqueness follows.

\subsubsection*{Endpoint continuity in time} To obtain the endpoint continuity, we regularize the initial data using the periodic heat kernel at time $t=\epsilon$ for certain $\epsilon>0$ fixed. The corresponding solution $(u^\epsilon,v^\epsilon)$ exists in 
$$
(u^\epsilon,v^\epsilon)\in L^\infty(0,T;H^4\times H^5)\cap C([0,T],H^3\times H^4),
$$ up to certain positive time that depends on the $C^1\times C^2$ semi-norms of the solutions. Thus, they are uniformly bounded in terms of the $H^3\times H^4$ size of the initial data. As a consequence, both solutions $(u,v)$ and $(u^\epsilon,v^\epsilon)$ are defined in the same lifespan. We consider now
$$
\sup_{t}(\|u(t)-u^\epsilon (t)\|_{H^3}+\|v(t)-v^\epsilon (t)\|_{H^4}).
$$
Performing the same kind of energy estimates as before (and the same as in the uniqueness step), we find that
$$
\sup_{t}(\|u(t)-u^\epsilon (t)\|_{H^3}+\|v(t)-v^\epsilon (t)\|_{H^4})\leq C\epsilon,
$$
from where the uniform convergence is obtained. The continuity follows.

\subsection{Proof of Theorem \ref{teobkdv2}}\label{proof:section:wp2}
As in the proof of Theorem \ref{teobkdv1}, let us focus on deriving the a priori energy estimates. Recall that the space $\mathbb{A}_{\mu(t)}^{s}$ contains analytic functions and that it enjoys the Banach Algebra property, namely
\[ \norm{fg}_{\mathbb{A}_{\mu(t)}^{s}}\leq \norm{f}_{\mathbb{A}_{\mu(t)}^{s}}\norm{g}_{\mathbb{A}_{\mu(t)}^{s}}.\]
We compute 
\begin{align*}
 \frac{d}{dt}\|u\|_{\mathbb{A}_{\mu(t)}}&=\frac{d}{dt} \displaystyle\sum_{k\in\mathbb{Z}} e^{\mu(t) |k|}|\widehat{u}(k,t)|= \displaystyle\sum_{k\in\mathbb{Z}} |k| \mu'(t)e^{\mu(t) |k|}|\widehat{u}(k,t)| +\displaystyle\sum_{k\in\mathbb{Z}} e^{\mu(t) |k|}\mathfrak{R}\bigg(\widehat{\partial_{t} u }(k,t)\frac{\overline{\widehat{u}(k,t)}}{|\widehat{u}(k,t)|} \bigg) \\
 &\leq \mu'(t) \displaystyle\sum_{k\in\mathbb{Z}} |k| e^{\mu(t) |k|}|\widehat{u}(k,t)| + \norm{u_{t}}_{\mathbb{A}_{\mu(t)}}
 \end{align*}
Therefore, taking $\mu(t)=\mu_{0}-Ct$ for a certain $C$ to be chosen and using equation \eqref{diss:BKdV:u} we find that
\begin{align}
 \frac{d}{dt}\|u\|_{\mathbb{A}_{\mu(t)}}&\leq -C\displaystyle\sum_{k\in\mathbb{Z}} |k| e^{\mu(t) |k|}|\widehat{u}(k,t)| + \norm{-3u\partial_x u-6v^2\partial_x v}_{\mathbb{A}_{\mu(t)}}.
\end{align}
Using the Banach Algebra property we observe that
\begin{align}
 \frac{d}{dt}\|u\|_{\mathbb{A}_{\mu(t)}}&\leq -C\norm{u}_{\mathbb{\dot{A}}^{1}_{\mu(t)}}+ 3 \norm{u}_{\mathbb{A}_{\mu(t)}}\norm{u}_{\mathbb{\dot{A}}^{1}_{\mu(t)}}
 +6\norm{v^2}_{\mathbb{A}_{\mu(t)}}\norm{v}_{\mathbb{\dot{A}}^{1}_{\mu(t)}}.
\end{align}
Repeating the same computations for the evolution of $v$, one can show that
\begin{align}
 \frac{d}{dt}\|v\|_{\mathbb{A}_{\mu(t)}}&\leq -C\norm{v}_{\mathbb{\dot{A}}^{1}_{\mu(t)}}+ 6 \norm{v}_{\mathbb{A}_{\mu(t)}}\norm{v}_{\mathbb{\dot{A}}^{1}_{\mu(t)}}
+\|u\|_{\mathbb{A}_{\mu(t)}}\|v\|_{\mathbb{\dot{A}}^{1}_{\mu(t)}}+\|v\|_{\mathbb{A}_{\mu(t)}}\|u\|_{\mathbb{\dot{A}}^{1}_{\mu(t)}}
\end{align}
Collecting both estimates and taking
$$
C=2(3\|u_0\|_{\mathbb{A}_{\mu_0}}+6\|v_0\|_{\mathbb{A}_{\mu_0}}^2+6\|v_0\|_{\mathbb{A}_{\mu_0}}+\|u_0\|_{\mathbb{A}_{\mu_0}}+\|v_0\|_{\mathbb{A}_{\mu_0}})
$$
we find that
$$
\frac{d}{dt}(\|u\|_{\mathbb{A}_{\mu(t)}}+\|v\|_{\mathbb{A}_{\mu(t)}})\leq 0.
$$
for certain uniform time $T=T(u_0,v_0,\mu_0)$. After this, the existence of solution follows from a standard mollifier
approach together with Picard’s theorem as shown in Theorem \ref{teobkdv1}. Similarly, the uniqueness follows from a standard contradiction argument. This concludes the result.

\section{Proof of the Theorem \ref{teomrs}}\label{proof:section:finite}
We argue by contradiction. Let us assume that the solution of \eqref{MRSLocal} remains smooth globally in time for initial data $u_0$ and $v_0$ satisfying the hypothesis stated in Theorem \ref{teomrs}. First notice that the odd symmetry is preserved for smooth solutions of the evolution equation \eqref{MRSLocal}. Define 
$$
U^{(\ell)}(t)=\partial_x^\ell u(x,t)\bigg{|}_{x=0},
$$
and
$$
V^{(\ell)}(t)=\partial_x^\ell v(x,t)\bigg{|}_{x=0}.
$$
Then, we compute that the solution of \eqref{MRSLocal} satisfy
\begin{align*}
\frac{d}{dt}U^{(1)}&=-(U^{(1)})^2+ V^{(1)},\\
\frac{d}{dt}V^{(1)}&=-(V^{(1)})^2-U^{(1)}.
\end{align*}
To show the finite time blow up of this system we define the functions
$$
F(t)=U^{(1)}-V^{(1)},
$$
$$
G(t)=U^{(1)}+V^{(1)}.
$$
Then we obtain that
\begin{align*}
\frac{d}{dt}F&=-FG+G,\\
\frac{d}{dt}G&=-G^2+2((F+G)/2)((G-F)/2)-F.
\end{align*}
As a consequence, we find that if $F(0)=1$, then $F(t)=1$ and the equation for $G$ becomes
$$
\frac{d}{dt}G=-G^2+(1+G)(G-1)/2-1=-G^2/2-3/2.
$$
From here we find that if $G(0)<-1$ the solution blows up in finite time and this is a contradiction with our initial hypotheses of global solvability.

\section*{Acknowledgement} 
D.A-O is supported by the fellowship of the Santander-ULL program.  D.A-O and R. G-B are also supported by the project “An\'alisis Matem\'atico Aplicado y Ecuaciones Diferenciales” Grant PID2022-141187NB-I00 and acronym "AMAED". This publication is part of the project PID2022-141187NB-I00 funded by MICIU/AEI /10.13039/501100011033 and by FEDER, UE/AEI /10.13039/501100011033 / FEDER, UE

 \end{document}